\documentclass[12pt, reqno]{amsart}
\allowdisplaybreaks
\begin{document}
\setcounter{page}{1}

\title[\hfilneg \hfil $\psi$ Hilfer fractional  ]
{Properties of  Some $\psi$-Hilfer fractional Integrodifferential equations}

\author[Deepak Pachpatte \hfil \hfilneg]
{Deepak B. Pachpatte}

\address{Deepak B. Pachpatte \newline
 Dept. of Mathematics,
 Dr. B. A. M. University, Aurangabad,
 Maharashtra 431004, India}
\email{pachpatte@gmail.com}

\subjclass[2010]{26A33, 34A08, 34A12, 34A40}

\keywords{ inequality, $\psi$ Hilfer fractional,estimate of solution, continuous dependence }

\begin{abstract}
 In this paper we study some properties of  $\psi$-Hilfer fractional integrodifferential equations.  We obtain the existence and uniqueness and other properties such as continuous dependence of solution. The tools used for obtaining our result is Banach fixed point theorem.  Certain  inequalities are also used obtain the estimates on the solution of the equation.
\end{abstract}
\maketitle
\section{Introduction}

   In last few decade abundance of applications are stimulating rapidly with the development of fractional differential equations. An excellent account of information on recent developments and its application in control theory and life sciences can be found in some recent books See \cite{Bal,Pet}
During last few years due to important applications in various fields fractional calculus has gained lot of attention from researchers. Recently many authors have studied the fractional integrodifferential equations using various definitions of fractional operators \cite{Kar}. The authors in \cite{Abb, Alm1, Fur,Sub,Tat,Tha, Wan, Zha} have studied various properties of fractional integrodifferential equations using various fractional operators.
  In \cite{Har1,Van2, Van3, Van5} various results are obtained for $\psi$-Hilfer fractional differential equations. Stability study of various $\psi$-Hilfer fractional differential and partial differential equations is done in \cite{Har2, Van3, Van4}.
  \paragraph{} Motivated by diverse applications and to widen the scope of such equations in this paper we  consider the $\psi$-Hilfer fractional integrodifferential equation of the form 
	\[
{}^HD_{a + }^{\alpha ,\beta ;\psi } \mathfrak{z}(x) = f\left( {x,\mathfrak{z}(x),\int\limits_a^x {w\left( {x,t,\mathfrak{z}(t)} \right)dt} } \right),
\tag{1.1}\]
\[
  I_{a + }^{1 - \gamma ,\psi } \mathfrak{z}\left( a \right) = \mathfrak{z}_a, \,\,
\tag{1.2}
\]
where
$
\gamma  = \alpha  + \beta \left( {1 - \alpha } \right)
$,      
 ${}^HD_{a + }^{\alpha ,\beta ;\psi } \mathfrak{z}(x) $ is the $\psi$-Hilfer fractional derivative  of order $\alpha$ and type $\beta$, $f:[a,b) \times \mathbb{R} \times \mathbb{R} \to \mathbb{R} $, $w:[a,b) \times [a,b) \times \mathbb{R} \to \mathbb{R} $, $\mathfrak{z}:C_{1 - \gamma } \left[ {a,b} \right]$ and $x \in [a,b]$.

	\section{Preliminaries}
	Now in this section we give the basic definitions and lemmas required in further subsequent discussions.
	Now similar to as given in \cite{Kil1, Van2,Van1, Van3}  $[a,b]$ with $(0<a<b<\infty )$ be a finite interval in $\mathbb{R}^+$ and $\mathcal{C}[a,b], A\mathcal{C}^{n}[a,b],\mathcal{C}^{n}[a,b]$ be the spaces of continuous function, n-times absolutely continuous and n-times continuously differentiable functions on $[a,b]$ respectively.
	In space of continuous function $\mathfrak{u}$ on $[a,b]$ is defined by
	\[
\left\| \mathfrak{u} \right\|_{\mathcal{C}\left[ {a,b} \right]}  = \mathop {\max }\limits_{x \in [a,b]} \left| {\mathfrak{u}\left( x \right)} \right|.
\]
The n-times absolutely continuous function $\mathfrak{u}$ on $[a,b]$ is defined as 
 \[
A\mathcal{C}^n \left[ {a,b} \right] = \left\{ {\mathfrak{u}:\left[ {a,b} \right] \to \mathbb{R};\mathfrak{u}^{\left( {n - 1} \right)}  \in A\mathcal{C}\left[ {a,b} \right]} \right\}.
\]
The weighted space $\mathcal{C}_{\gamma ;\psi} [a,b] $ of function $u$ on (a,b] is defined by
	\[
\mathcal{C}_{\gamma ;\psi} [a,b]  = \left\{ {\mathfrak{u}:(a,b] \to \mathbb{R},\left( {\psi \left( x \right) - \psi \left( a \right)} \right)^\gamma,  \mathfrak{u}\left( x \right) \in \mathcal{C}\left[ {a,b} \right]} \right\},0 \le \gamma  < 1,
\]
with norm
\[
\left\| \mathfrak{u}\right\|_{\mathcal{C}_{\gamma ;\psi} [a,b] }  = \left\| {\left( {\psi \left( x \right) - \psi \left( a \right)} \right)^\gamma  \mathfrak{u}(x)} \right\|_{\mathcal{C}[a,b]}  = \mathop {\max }\limits_{x \in \left[ {a,b} \right]} \left| {\left( {\psi \left( x \right) - \psi \left( a \right)} \right)^\gamma  \mathfrak{u}\left( x \right)} \right|.
\]
The weighted space $\mathcal{C}_{\gamma ;\psi}^n [a,b] $ of a function $u$ on $(a,b]$ is defined by
\[
\mathcal{C}_{\gamma ;\psi }^n \left[ {a,b} \right] = \left\{ {\mathfrak{u}:(a,b] \to \mathbb{R},\,\,\,\mathfrak{u}\left( x \right) \in \mathcal{C}^{n - 1} \left[ {a,b} \right],\,\mathfrak{u}^{(n)} \left( x \right) \in \mathcal{C}_{\gamma ;\psi } [a,b]\,} \right\},\,\,0 \le \gamma  < 1,
\]
with the norm
\[
\left\| \mathfrak{u}\right\|_{\mathcal{C}_{\gamma ;\psi }^n [a,b]}  = \sum\limits_{k = 0}^{n - 1} {\left\| {c^{\left( k \right)} } \right\|} _{\mathcal{C}[a,b]}  + \left\| {c^{\left( n \right)} } \right\|_{\mathcal{C}_{\gamma ;\psi [a,b]} } .
\]
The weighted space $\mathcal{C}_{1 - \gamma ;\psi } [a,b]$ of a function $\mathfrak{u}$ on [a,b]
\[
\mathcal{C}_{1 - \gamma ;\psi } [a,b] = \left\{ {\mathfrak{u}:\left( {a,b} \right] \to \mathbb{R},\left( {\psi \left( x \right) - \psi \left( a \right)} \right)^{1 - \gamma } \mathfrak{u}\left( t \right) \in \mathcal{C}[a,b]} \right\} ,\,\,0 \le \gamma  < 1,
\]
and
\begin{align*}
\left\| \mathfrak{u}\right\|_{\mathcal{C}_{1 - \gamma ;\psi } [a,b]} 
& = \left\| {\left( {\psi \left( x \right) - \psi \left( a \right)} \right)^{1 - \gamma } \mathfrak{z}\left( t \right)} \right\|_{\mathcal{C}[a,b]}  \\
& = \mathop {\max }\limits_{x \in (a,b]} \left| {\left( {\psi \left( x \right) - \psi \left( a \right)} \right)^{1 - \gamma } \mathfrak{u}\left( t \right)} \right|.
\end{align*}
For $n=0$, we get $\mathcal{C}_\gamma ^0 \left[ {a,b} \right] = \mathcal{C}_\gamma  \left[ {a,b} \right]. $ \\
The weighted space $\mathcal{C}_{\gamma ;\psi }^{\alpha  , \beta } [a,b]$ is defined by
\[
\mathcal{C}_{\gamma ;\psi }^{\alpha  , \beta } [a,b] = \left\{ {\mathfrak{u}\in \mathcal{C}_{\gamma ;\psi } [a,b]:{}^HD_{a + }^{\alpha ,\beta ;\psi } \mathfrak{u}\in \mathcal{C}_{\gamma ;\psi } [a,b]\,} \right\},\,\,\gamma  = \alpha  + \beta \left( {1 - \alpha } \right).
\]
In \cite{Kil1,Sam1} the authors have defined the fractional integrals and fractional derivative of a function with respect to another function as follows:
\paragraph{\textbf{Definition 2.1 }} \cite{Alm1,Kil1}. Let $I=[a,b]$ be an interval, $\alpha >0$, $f$ is an integrable function defined on $I$ and $\psi \in \mathcal{C}^1(I)$ an increasing function such that $\psi '\left( x \right) \ne 0$ for all $x \in I$ then fractional derivative and integral of $f$ is given by
\[
I_{a + }^{\alpha ,\psi } f(x) = \frac{1}{{\Gamma \left( \alpha  \right)}}\int\limits_a^x {\psi '\left( t \right)\left( {\psi \left( x \right) - \psi \left( t \right)} \right)^{\alpha  - 1} f\left( t \right)dt} 
\]
 and
\begin{align*}
 D_{a + }^{\alpha ,\psi } f\left( x \right)
&= \left( {\frac{1}{{\psi '\left( x \right)}}\frac{d}{{dx}}} \right)^n I_{a + }^{n - \alpha ,\psi } f\left( x \right) \\ 
&= \frac{1}{{\Gamma \left( {n - \alpha } \right)}}\left( {\frac{1}{{\psi '\left( x \right)}}\frac{d}{{dx}}} \right)^n \int\limits_a^x {\psi '\left( t \right)\left( {\psi \left( x \right) - \psi \left( t \right)} \right)^{n - \alpha  - 1} f\left( t \right)dt}, 
\end{align*}
respectively. Similarly right fractional integral and right fractional derivative are given by
\[
I_{b - }^{\alpha ,\psi } f(x) = \frac{1}{{\Gamma \left( \alpha  \right)}}\int\limits_a^x {\psi '\left( t \right)\left( {\psi \left( t \right) - \psi \left( x \right)} \right)^{\alpha  - 1} f\left( t \right)dt} 
\]
 and
\begin{align*}
 D_{b- }^{\alpha ,\psi } f\left( x \right)
&= \left( -{\frac{1}{{\psi '\left( x \right)}}\frac{d}{{dx}}} \right)^n I_{b - }^{n - \alpha ,\psi } f\left( x \right) \\ 
&= \frac{1}{{\Gamma \left( {n - \alpha } \right)}}\left( {\frac{1}{{\psi '\left( x \right)}}\frac{d}{{dx}}} \right)^n \int\limits_a^x {\psi '\left( t \right)\left( {\psi \left( t \right) - \psi \left( x \right)} \right)^{n - \alpha  - 1} f\left( t \right)dt}. 
\end{align*}
\paragraph{}In \cite{Van1} authors have defined $\psi$-Hilfer fractional derivative as follows:
\paragraph{\textbf{Definition 2.2 }}
Let $n - 1 < \alpha  < n $ with $n \in N$, $[a,b]$ is the interval such that $ - \infty  \le a < b \le \infty  $ and $y,\psi  \in \mathcal{C}^n \left( {\left[ {a,b} \right],\mathbb{R}} \right)$ two functions such that $\psi$ is increasing and $\psi '\left( x \right) \ne 0$, for all $x \in \left[ {a,b} \right] $. The $\psi$-Hilfer derivative (left sided and right sided)  ${}^HD_{a + }^{\alpha ,\beta ;\psi } \left( . \right)$ of a function of order $\alpha$ and type $0 \le \beta  \le 1 $ are defined by
\[
{}^HD_{a + }^{\alpha ,\beta ;\psi } y(x) = I_{a + }^{\beta \left( {n - \alpha } \right);\psi } \left( {\frac{1}{{\psi '\left( x \right)}}\frac{d}{{dx}}} \right)^n I_{a + }^{\left( {1 - \beta } \right)\left( {1 - \alpha } \right):\psi } y\left( x \right)
\]
and
\[
{}^HD_{b - }^{\alpha ,\beta ;\psi } y(x) = I_{b - }^{\beta \left( {n - \alpha } \right);\psi } \left( { - \frac{1}{{\psi '\left( x \right)}}\frac{d}{{dx}}} \right)^n I_{b - }^{\left( {1 - \beta } \right)\left( {1 - \alpha } \right):\psi } y\left( x \right).
\]
\paragraph{}In particular when  $0 < \alpha  < 1$ and $0 \le \beta  \le 1$ we get 
\[
{}^HD_{a + }^{\alpha ,\beta ;\psi } y(x) = \frac{1}{{\Gamma \left( {\gamma  - \alpha } \right)}}\int\limits_a^x {\left( {\psi (x) - \psi (t)} \right)^{\gamma  - \alpha  - 1} D_{a + }^{\gamma ;\psi } y\left( t \right)dt}; 
\]
with $\gamma  = \alpha  + \beta \left( {1 - \alpha } \right)$ and $D_{a + }^{\gamma ;\psi } (.)$ is $\psi$-Riemann-liouville fractional derivative.\\
Now we give some theorems proved in \cite{Van1} which will be required in proving our results.
\paragraph{\textbf{Lemma 2.1}} [\cite{Van1} Lemma 4, p 80]
Let $n - 1 \le \gamma  < n$ and $f \in \mathcal{C}_\gamma  \left[ {a,b} \right]$. Then
\[
I_{a + }^{\alpha ;\psi } f\left( a \right) = \mathop {\lim }\limits_{x \to a + } I_{a + }^{\alpha ;\psi } f\left( x \right) = 0,\,\,\,n - 1 \le \gamma  < \alpha. 
\].
\paragraph{\textbf{Lemma 2.2}} [\cite{Van1} Theorem 7, p 80]
Let $f \in \mathcal{C}^1 \left[ {a,b} \right]$, $\alpha  > 0$ and $0 \le \beta  \le 1$, we have \\
${}^HD_{a + }^{\alpha ,\beta ;\psi } I_{a + }^{\alpha ;\psi } f\left( x \right) = f\left( x \right)$ and ${}^HD_{b - }^{\alpha ,\beta ;\psi } I_{b - }^{\alpha ;\psi } f\left( x \right) = f\left( x \right)$.\\
 \paragraph{} In \cite{Bgp1} Pachpatte has proved  Gronwall-Bellman inequality  is as follows:
\paragraph{\textbf{Lemma 2.3}}
Let $\mathfrak{u}(t), f(t)$ and $g(t)$ be real valued nonnegative continuous functions defined on $[0,\infty)$, for which the inequality 
\[
\mathfrak{u}\left( t \right) \le \mathfrak{u}_0  + \int\limits_0^t {f(s)\mathfrak{u}(s)ds + } \int\limits_0^t {f(s)\left( {\int\limits_0^s {g\left( \tau  \right)\mathfrak{u}\left( \tau  \right)d\tau } } \right)}, 
\]
for $t \in [0,\infty)$, holds, where $\mathfrak{u}_0$ is a nonnegative constant. Then
\[
\mathfrak{u}\left( t \right) \le \mathfrak{u}_0 \left( {1 + \int\limits_0^t {f(s)\exp \left( {\int\limits_0^s {\left( {f(\tau ) + g(\tau )} \right)d\tau } } \right)ds} } \right)
\]
for $t \in [0,\infty)$.

Now we define  what is mean by a solution of the problem $(1.1)$.
\paragraph{\textbf{Definition 2.3 }} A function $\mathfrak{z} \in \mathcal{C}[a,b]$ is said to be a solution of $(1.1)$ if $\mathfrak{z}$ satisfies equation $(1.1)$ and the conditions in $(1.2)$.
\paragraph{\textbf{Lemma 2.4}}
	A function $\mathfrak{z} \in \mathcal{C}[a,b]$ is a solution of $(1.1)-(1.2)$ if and only if $\mathfrak{z}$ satisfies
\begin{align*}
\mathfrak{z}\left( x \right)
&= \frac{{\mathfrak{z}_a }}{{\Gamma \left( \gamma  \right)}}\left( {\psi (x) - \psi (a)} \right)^{\gamma  - 1}  \\ 
&+ \frac{1}{{\Gamma \left( \alpha  \right)}}\int\limits_a^x {\psi '(t)\left( {\psi (x) - \psi (t)} \right)^{\alpha  - 1} f\left( {t,\mathfrak{z}(t),\int\limits_a^t {w\left( {t,\tau ,\mathfrak{z}\left( \tau  \right)} \right)d\tau } } \right)} dt 
\tag{2.1}
\end{align*}
that is $(1.1)$ and $(2.1)$
 are equivalent. 
\paragraph{\textbf{Proof.}} Let $\mathfrak{z}(x)$ be a solution of $(1.1)-(1.2)$. Then from the definition of $\psi$-Hilfer fractional derivative ${}^HD_{a + }^{\alpha ,\beta ;\psi } $ we have
\[
I_{a + }^{\alpha ;\psi } {}^HD_{a + }^{\alpha ,\beta ;\psi } \mathfrak{z}(x) = \mathfrak{z}(x) - \frac{{\left( {\psi (x) - \psi (a)} \right)^{\gamma  - 1} }}{{\Gamma \left( \gamma  \right)}}I_{a + }^{\left( {1 - \beta } \right)\left( {1 - \alpha } \right);\psi } \mathfrak{z}(a),
\tag{2.2}\]
for $x \in [a,b]$. Thus using Lemma  $(2.2)$ in  $(2.2)$ we have 
\begin{align*}
\mathfrak{z}\left( x \right)
&= \frac{{\mathfrak{z}_a }}{{\Gamma \left( \gamma  \right)}}\left( {\psi (x) - \psi (a)} \right)^{\gamma  - 1}  \\ 
&+ \frac{1}{{\Gamma \left( \alpha  \right)}}\int\limits_a^x {\psi '(t)\left( {\psi (x) - \psi (t)} \right)^{\alpha  - 1} f\left( {t,\mathfrak{z}(t),\int\limits_a^t {w\left( {t,\tau ,\mathfrak{z}\left( \tau  \right)} \right)d\tau } } \right)} dt 
\end{align*}
which is equation $(2.1)$.
\paragraph{} Apply the operator $D_{a + }^{\alpha ,\beta ;\psi } $ on both sides of $(2.1)$ we have
\begin{align*}
{}^HD_{a + }^{\alpha ,\beta ;\psi } \mathfrak{z}(x)
&= {}^HD_{a + }^{\alpha ,\beta ;\psi } \left[ {\frac{{\mathfrak{z}_a }}{{\Gamma \left( \gamma  \right)}}\left( {\psi (x) - \psi (a)} \right)^{\gamma  - 1} } \right. \\ 
& \left. { + \frac{1}{{\Gamma \left( \alpha  \right)}}\int\limits_a^x {\psi '(t)\left( {\psi (x) - \psi (t)} \right)^{\alpha  - 1} f\left( {t,\mathfrak{z}(t),\int\limits_a^t {w\left( {t,\tau ,\mathfrak{z}\left( \tau  \right)} \right)d\tau } } \right)} dt} \right].
\end{align*}
Using Lemma $(2.1)$ and Lemma $(2.2)$ as $x \rightarrow a$ we have
\[
{}^HD_{a + }^{\alpha ,\beta ;\psi } \mathfrak{z}(x) = f\left( {x,\mathfrak{z}\left( x \right),\int\limits_a^x {w\left( {x,t,\mathfrak{z}(t)} \right)dt} } \right),
\]
which is $(1.1)$.
	\section{Existence and Uniqueness }
	Now in our next theorem we prove the existence and uniqueness of solution of equation $(1.1)-(1.2)$.
	\paragraph{\textbf{Theorem 3.1.}} Let $0 < \alpha  < 1$, $0 \le \beta  \le 1$ and $\gamma  = \alpha  + \beta  - \alpha \beta$. Suppose the functions $f,w$ in $(1.1)$ be continuous and satisfy the conditions
	\[
\left| {f\left( {s,t,g} \right) - f\left( {s,\overline t ,\overline g } \right)} \right| \le Q_1 \left[ {\left| {t - \overline t } \right| + \left| {g - \overline g } \right|} \right],
\tag{3.1}\]
\[
\left| {w\left( {s,g,\eta } \right) - w\left( {s,g,\overline \eta  } \right)} \right| \le Q_2 \left| {\eta  - \overline \eta  } \right|
\tag{3.2}\]
and let 
\begin{align*}
p 
& = \mathop {\max }\limits_{x \in [a,b]} \left( {\psi \left( x \right) - \psi \left( a \right)} \right)^{1 - \gamma } \left| {\frac{{\mathfrak{z}_a }}{{\Gamma \left( \alpha  \right)}}\left( {\psi \left( x \right) - \psi \left( a \right)} \right)^{\gamma  - 1} } \right. \\
&\left. { + I_{a + }^{\alpha ;\psi } f\left( {x,0,\int\limits_a^x {w\left( {x,t,0} \right)dt} } \right)} \right|  
 < \infty.
\tag{3.3}
\end{align*}
If
\begin{align*}
	q
	& = Q_1 \frac{{\Gamma \left( \gamma  \right)}}{{\Gamma \left( {\alpha  + \gamma } \right)}}\left( {\psi \left( b \right) - \psi \left( a \right)} \right)^\alpha    \\ 
	&+ Q_1 Q_2 \frac{{\Gamma \left( \gamma  \right)}}{{\Gamma \left( {\alpha  + \gamma } \right)}}\left( {\psi \left( b \right) - \psi \left( a \right)} \right)^{\alpha  + 1}  < 1,	
	\tag{3.4}
	\end{align*}
	then equation $(1.1)$ has a unique solution $\mathfrak{z} \in \mathcal{C}[a,b]$.
\paragraph{\textbf{Proof.}} Define the operator $\mathfrak{F}$ by
\begin{align*}
\left( { \mathfrak{F}\mathfrak{z}} \right)\left( x \right)
& = \frac{{\mathfrak{z}_a }}{{\Gamma \left( \gamma  \right)}}\left( {\psi \left( x \right) - \psi \left( a \right)} \right)^{\gamma  - 1}  + I_{a + }^{\alpha ;\psi } f\left( {t,\mathfrak{z}\left( t \right),\int\limits_a^t {w\left( {t,\tau ,\mathfrak{z}\left( \tau  \right)} \right)d\tau } } \right) \\
&- I_{a + }^{\alpha ;\psi } f\left( {t,0,\int\limits_a^t {w\left( {t,\tau ,0} \right)d\tau } } \right) + I_{a + }^{\alpha ;\psi } f\left( {t,0,\int\limits_a^t {w\left( {t,\tau ,0} \right)d\tau } } \right), 
\tag{3.5}
\end{align*}
for $x \in [a,b]$.\\ Now we prove that $\mathfrak{F}$ is a contraction map. From equation $(3.5)$ and hypotheses we have
\begin{align*}
\left\| {\mathfrak{F}\mathfrak{z}} \right\|_{\mathcal{C}_{1 - \gamma ,\psi } [a,b]}
&= \mathop {\max }\limits_{x \in (a,b]} \left( {\psi \left( x \right) - \psi \left( a \right)} \right)^{1 - \gamma } \left| {\mathfrak{F}\mathfrak{z}} \right| \\
& = \mathop {\max }\limits_{x \in (a,b]} \left( {\psi \left( x \right) - \psi \left( a \right)} \right)^{1 - \gamma }  \\
& \times \left| {\frac{{\mathfrak{z}_a }}{{\Gamma \left( \gamma  \right)}}\left( {\psi \left( x \right) - \psi \left( a \right)} \right)^{\gamma  - 1}  + I_{a + }^{\alpha ;\psi } f\left( {t,\mathfrak{z}\left( t \right),\int\limits_a^t {w\left( {t,\tau ,\mathfrak{z}\left( \tau  \right)} \right)d\tau } } \right)} \right. \\
&\left. { - I_{a + }^{\alpha ;\psi } f\left( {t,0,\int\limits_a^t {w\left( {t,\tau ,0} \right)d\tau } } \right) + I_{a + }^{\alpha ;\psi } f\left( {t,0,\int\limits_a^t {w\left( {t,\tau ,0} \right)d\tau } } \right)} \right| \\ 
& \le \mathop {\max }\limits_{x \in (a,b]} \left( {\psi \left( x \right) - \psi \left( a \right)} \right)^{1 - \gamma }  \\
& \times \left| {\,\,\frac{{\mathfrak{z}_a }}{{\Gamma \left( \gamma  \right)}}\left( {\psi \left( x \right) - \psi \left( a \right)} \right)^{\gamma  - 1}  + I_{a + }^{\alpha ;\psi } f\left( {t,0,\int\limits_a^t {w\left( {t,\tau ,0} \right)d\tau } } \right)\,} \right| \\ 
& + \mathop {\max }\limits_{x \in (a,b]} \left( {\psi \left( x \right) - \psi \left( a \right)} \right)^{1 - \gamma }  \\ 
&\times I_{a + }^{\alpha ;\psi } \left| {f\left( {t,\mathfrak{z}\left( t \right),\int\limits_a^t {w\left( {t,\tau ,\mathfrak{z}\left( \tau  \right)} \right)d\tau } } \right) - f\left( {t,0,\int\limits_a^t {w\left( {t,\tau ,0} \right)d\tau } } \right)} \right| \\
& \le p + \mathop {\max }\limits_{x \in (a,b]} \left( {\psi \left( x \right) - \psi \left( a \right)} \right)^{1 - \gamma }  \\
&  \times I_{a + }^{\alpha ;\psi } Q_1 \left[ {\left| {\mathfrak{z}\left( t \right)} \right| + \int\limits_a^t {\left| {w\left( {t,\tau ,\mathfrak{z}\left( \tau  \right)} \right) - w\left( {t,\tau ,0} \right)} \right|d\tau } } \right] \\ 
&\le p + \mathop {\max }\limits_{x \in (a,b]} \left( {\psi \left( x \right) - \psi \left( a \right)} \right)^{1 - \gamma } I_{a + }^{\alpha ;\psi } Q_1 \left[ {\left| {\mathfrak{z}\left( t \right)} \right| + Q_2 \int\limits_a^t {\left| {\mathfrak{z}\left( \tau  \right)} \right|d\tau } } \right]\\
&\le p + \mathop {\max }\limits_{x \in (a,b]} \left( {\psi \left( x \right) - \psi \left( a \right)} \right)^{1 - \gamma }  \\ 
& \times I_{a + }^{\alpha ;\psi } Q_1 \left[ {\frac{{\left| {\mathfrak{z}\left( t \right)} \right|}}{{\left( {\psi \left( x \right) - \psi \left( a \right)} \right)^{1 - \gamma } }}\left( {\psi \left( x \right) - \psi \left( a \right)} \right)^{1 - \gamma } } \right. \\ 
&\left. { + Q_2 \int\limits_a^t {\frac{{\left| {\mathfrak{z}\left( t \right)} \right|}}{{\left( {\psi \left( x \right) - \psi \left( a \right)} \right)^{1 - \gamma } }}\left( {\psi \left( x \right) - \psi \left( a \right)} \right)^{1 - \gamma }d\tau } } \right] \\ 
&\le p + \left\| \mathfrak{z} \right\|_{\mathcal{C}_{1 - \gamma ;\psi } [a,b]} \mathop {\max }\limits_{x \in (a,b]} \left( {\psi \left( x \right) - \psi \left( a \right)} \right)^{1 - \gamma }  \\ 
&  \times Q_1 I_{a + }^{\alpha ;\psi
} \left[ {\left( {\psi \left( x \right) - \psi \left( a \right)} \right)^{\gamma  - 1}  + Q_2 \int\limits_a^t {\left( {\psi \left( x \right) - \psi \left( a \right)} \right)^{\gamma  - 1} d\tau } } \right] \\ 
&= p + \left\| \mathfrak{z} \right\|_{\mathcal{C}_{1 - \gamma ;\psi } [a,b]} \mathop {\max }\limits_{x \in (a,b]} \left( {\psi \left( x \right) - \psi \left( a \right)} \right)^{1 - \gamma }  \\
& \times Q_1 I_{a + }^{\alpha ;\psi } \left[ {\left( {\psi \left( x \right) - \psi \left( a \right)} \right)^{\gamma  - 1}  + \frac{{Q_2 }}{\gamma }\int\limits_a^t {\left( {\psi \left( x \right) - \psi \left( a \right)} \right)^\gamma  d\tau } } \right] \\ 
& = p + \left\| \mathfrak{z} \right\|_{\mathcal{C}_{1 - \gamma ;\psi } [a,b]} \mathop {\max }\limits_{x \in (a,b]} \left( {\psi \left( x \right) - \psi \left( a \right)} \right)^{1 - \gamma }  \\
&\times \left[ {Q_1 \frac{{\Gamma \left( \gamma  \right)}}{{\Gamma \left( {\alpha  + \gamma } \right)}}\left( {\psi \left( x \right) - \psi \left( a \right)} \right)^{\alpha  + \gamma  - 1} } \right. \\ 
&\left. { + Q_1 Q_2 \frac{{\Gamma \left( \gamma  \right)}}{{\Gamma \left( {\alpha  + \gamma  + 1} \right)}}\left( {\psi \left( x \right) - \psi \left( a \right)} \right)^{\alpha  + \gamma } } \right] \\ 
&= p + \left\| \mathfrak{z} \right\|_{\mathcal{C}_{1 - \gamma ;\psi } [a,b]} \mathop {\max }\limits_{x \in (a,b]} \left[ {Q_1 \frac{{\Gamma \left( \gamma  \right)}}{{\Gamma \left( {\alpha  + \gamma } \right)}}\left( {\psi \left( x \right) - \psi \left( a \right)} \right)^\alpha  } \right. \\ 
&\left. { + Q_1 Q_2 \frac{{\Gamma \left( \gamma  \right)}}{{\Gamma \left( {\alpha  + \gamma  + 1} \right)}}\left( {\psi \left( x \right) - \psi \left( a \right)} \right)^{\alpha  + 1} } \right] \\ 
&\le p + \left\| \mathfrak{z} \right\|_{\mathcal{C}_{1 - \gamma ;\psi } [a,b]} \left[ {Q_1 \frac{{\Gamma \left( \gamma  \right)}}{{\Gamma \left( {\alpha  + \gamma } \right)}}\left( {\psi \left( b \right) - \psi \left( a \right)} \right)^\alpha  } \right. \\
&\left. { + Q_1 Q_2 \frac{{\Gamma \left( \gamma  \right)}}{{\Gamma \left( {\alpha  + \gamma  + 1} \right)}}\left( {\psi \left( b \right) - \psi \left( a \right)} \right)^{\alpha  + 1} } \right] \\
& = p + \left\| \mathfrak{z} \right\|_{\mathcal{C}_{1 - \gamma ;\psi } [a,b]} q \\
&< \infty. 
\tag{3.6}
\end{align*}
This proves that the operator $\mathfrak{F}$ maps $\mathcal{C}[a,b]$ into itself. 
\paragraph{}Now we verify that $\mathfrak{F}$ is a contraction map. 
Let $\mathfrak{z}_1,\mathfrak{z}_2 \in \mathcal{C}[a,b]$. From $(3.5)$ and hypotheses we get
\begin{align*}
&\left\| {\left( {\mathfrak{F}\mathfrak{z}_1 } \right)\left( x \right) - \left( {\mathfrak{F}\mathfrak{z}_2 } \right)\left( x \right)} \right\|_{\mathcal{C}_{1 - \gamma ,\psi } [a,b]} \\
&= \mathop {\max }\limits_{x \in \left[ {a,b} \right]} \left( {\psi \left( x \right) - \psi \left( a \right)} \right)^{1 - \gamma } \left| {\left( {\mathfrak{F}\mathfrak{z}_1 } \right)\left( x \right) - \left( {\mathfrak{F}\mathfrak{z}_2 } \right)\left( x \right)} \right| \\
& = \mathop {\max }\limits_{x \in \left[ {a,b} \right]} \left( {\psi \left( x \right) - \psi \left( a \right)} \right)^{1 - \gamma } \left| {I_a^{\alpha ;\psi } f\left( {t,\mathfrak{z}_1 \left( \tau  \right),\int\limits_a^t {w\left( {t,\tau ,\mathfrak{z}_1 \left( \tau  \right)} \right)d\tau } } \right)} \right. \\ 
&\left. { - I_a^{\alpha ;\psi } f\left( {t,\mathfrak{z}_2 \left( \tau  \right),\int\limits_a^t {w\left( {t,\tau ,\mathfrak{z}_2 \left( \tau  \right)} \right)d\tau } } \right)} \right| \\
 & \le \mathop {\max }\limits_{x \in \left[ {a,b} \right]} \left( {\psi \left( x \right) - \psi \left( a \right)} \right)^{1 - \gamma } Q_1 I_a^{\alpha ;\psi } \left[ {\left| {\mathfrak{z}_1 \left( t \right) - \mathfrak{z}_2 \left( t \right)} \right| + Q_2 \int\limits_a^t {\left| {\mathfrak{z}_1 \left( t \right) - \mathfrak{z}_2 \left( t \right)} \right|d\tau } } \right] \\
&= \mathop {\max }\limits_{x \in [a,b]} \left( {\psi \left( x \right) - \psi \left( a \right)} \right)^{1 - \gamma } Q_1 I_a^{\alpha ;\psi } \left[ {\frac{{\left| {\mathfrak{z}_1 \left( t \right) - \mathfrak{z}_2 \left( t \right)} \right|}}{{\left( {\psi \left( x \right) - \psi \left( a \right)} \right)^{1 - \gamma } }}} \right.\left( {\psi \left( x \right) - \psi \left( a \right)} \right)^{\gamma  - 1}  \\ 
& \left. { + Q_2 \int\limits_a^t {\frac{{\left| {\mathfrak{z}_1 \left( t \right) - \mathfrak{z}_2 \left( t \right)} \right|}}{{\left( {\psi \left( x \right) - \psi \left( a \right)} \right)^{1 - \gamma } }}} \left( {\psi \left( x \right) - \psi \left( a \right)} \right)^{\gamma  - 1} d\tau } \right] \\
&\le \mathop {\max }\limits_{x \in [a,b]} \left( {\psi \left( x \right) - \psi \left( a \right)} \right)^{1 - \gamma } Q_1 I_a^{\alpha ;\psi } \left[ {\left\| {\mathfrak{z}_1  - \mathfrak{z}_2 } \right\|_{\mathcal{C}_{1 - \gamma ,\psi \left[ {a,b} \right]} } } \right.\left( {\psi \left( x \right) - \psi \left( a \right)} \right)^{\gamma  - 1}  \\ 
&\left. { + Q_2 \int\limits_a^t {\left\| {\mathfrak{z}_1  - \mathfrak{z}_2 } \right\|_{\mathcal{C}_{1 - \gamma ,\psi \left[ {a,b} \right]} } \left( {\psi \left( x \right) - \psi \left( a \right)} \right)^{\gamma  - 1} d\tau } } \right] \\ 
& = \mathop {\max }\limits_{x \in [a,b]} \left( {\psi \left( x \right) - \psi \left( a \right)} \right)^{1 - \gamma } Q_1 \left\| {\mathfrak{z}_1  - \mathfrak{z}_2 } \right\|_{\mathcal{C}_{1 - \gamma ,\psi \left[ {a,b} \right]} }  \\
&\times I_a^{\alpha ;\psi } \left[ {\left( {\psi \left( x \right) - \psi \left( a \right)} \right)^{\gamma  - 1}  + Q_2 \int\limits_a^t {\left( {\psi \left( x \right) - \psi \left( a \right)} \right)^{\gamma  - 1} d\tau } } \right] \\ 
& = \left\| {\mathfrak{z}_1  - \mathfrak{z}_2 } \right\|_{\mathcal{C}_{1 - \gamma ,\psi \left[ {a,b} \right]} } \mathop {\max }\limits_{x \in [a,b]}  \\ 
& \left[ {Q_1 \frac{{\Gamma \left( \gamma  \right)}}{{\Gamma \left( {\alpha  + \gamma } \right)}}\left( {\psi \left( x \right) - \psi \left( a \right)} \right)^\alpha   + Q_1 Q_2 \frac{{\Gamma \left( \gamma  \right)}}{{\Gamma \left( {\alpha  + \gamma  - 1} \right)}}\left( {\psi \left( x \right) - \psi \left( a \right)} \right)^{\alpha  + 1} } \right] \\ 
& \le \left\| {\mathfrak{z}_1  - \mathfrak{z}_2 } \right\|_{\mathcal{C}_{1 - \gamma ,\psi \left[ {a,b} \right]} } \left[ {Q_1 \frac{{\Gamma \left( \gamma  \right)}}{{\Gamma \left( {\alpha  + \gamma } \right)}}\left( {\psi \left( b \right) - \psi \left( a \right)} \right)^\alpha  } \right. \\ 
&\left. { + Q_1 Q_2 \frac{{\Gamma \left( \gamma  \right)}}{{\Gamma \left( {\alpha  + \gamma  - 1} \right)}}\left( {\psi \left( b \right) - \psi \left( a \right)} \right)^{\alpha  + 1} } \right] \\ 
&= \left\| {\mathfrak{z}_1  - \mathfrak{z}_2 } \right\|_{\mathcal{C}_{1 - \gamma ,\psi \left[ {a,b} \right]} } q .
\tag{3.7}
\end{align*}
Since $q<1$ it follows from Banach fixed point theorem that $\mathfrak{F}$ has a unique fixed point in $\mathcal{C}[a,b]$. The fixed point of $\mathfrak{F}$ gives solution of equation $(1.1)$
\section{Estimate of Solution}
Many a times in order to study the various phenomena of the system,  obtaining estimates of the solutions is very important and is useful in the  analysis of the system. Now in our next theorem we obtain the estimates of solution of equation $(1.1)$.
\paragraph{\textbf{Theorem 4.1.}}
Suppose the functions $f,w$ be as in $(1.1)$ are continuous and satisfies the condition
\[
\left| {f\left( {x,t,g} \right) - f\left( {x,\overline t ,\overline g } \right)} \right| \le Q_3 \left( x \right)\left[ {\left| {t - \overline t } \right| + \left| {g - \overline g } \right|} \right],
\tag{4.1}\]
\[
\left| {w\left( {x,\upsilon ,\eta } \right) - w\left( {x,\upsilon ,\overline \eta  } \right)} \right| \le Q_4 \left( x \right)\left[ {\left| {\eta  - \overline \eta  } \right|} \right],
\tag{4.2}\]
where $Q_3 ,Q_4  \in \mathcal{C}\left( {[a,b]  ,\mathbb{R}_ +  } \right)$. Let 
\[
p_2  = \mathop {\max }\limits_{x \in [a,b]} \left| {\left| {\frac{{\mathfrak{z}_a }}{{\Gamma \left( \alpha  \right)}}} \right.\left( {\psi \left( x \right) - \psi \left( a \right)} \right)^{\gamma  - 1}  + I_{a + }^{\alpha ;\psi } f\left( {t,0,\int\limits_a^t {w\left( {t,\tau ,0} \right)d\tau } } \right)} \right|
\tag{4.3}\]
If $u(x), x\in [a,b]$ is any solution of $(1.1)$ then
\begin{align*}
\left| {\mathfrak{u}\left( x \right)} \right|
& \le p_2 \left[ {1 + \frac{1}{{\Gamma \left( \alpha  \right)}}} \right.\int\limits_a^x {\psi '\left( t \right)} \left( {\psi \left( x \right) - \psi \left( t \right)} \right)^{\alpha  - 1} Q_3 \left( t \right) \\
&\left. {\exp \left( {\int\limits_a^t {\left( {\frac{{\psi '\left( s \right)\left( {\psi \left( x \right) - \psi \left( s \right)} \right)^{\alpha  - 1} }}{{\Gamma \left( \alpha  \right)}}Q_3 \left( s \right) + Q_4 \left( s \right)} \right)ds} } \right)dt} \right]. \\
\tag{4.4}
\end{align*}

\paragraph{\textbf{Proof.}}
Since $\mathfrak{z}(x)$ is a solution of $(1.1)$ and from the 
\begin{align*}
 \left| {\mathfrak{z}\left( x \right)} \right| 
&\le \left| {\frac{{\mathfrak{z}_a }}{{\Gamma \left( \alpha  \right)}}} \right.\left( {\psi \left( x \right) - \psi \left( a \right)} \right)^{\gamma  - 1}  \\ 
&\left. { + \frac{1}{{\Gamma \left( \alpha  \right)}}\int\limits_a^x {\psi '\left( t \right)} \left( {\psi \left( x \right) - \psi \left( t \right)} \right)^{\alpha  - 1} f\left( {t,0,\int\limits_a^t {w\left( {t,\tau ,0} \right)d\tau } } \right)} \right| \\
& + \left| {\frac{1}{{\Gamma \left( \alpha  \right)}}\int\limits_a^x {\psi '\left( t \right)} \left( {\psi \left( x \right) - \psi \left( t \right)} \right)^{\alpha  - 1} f\left( {t,\mathfrak{z}(t),\int\limits_a^t {w\left( {t,\tau ,\mathfrak{z}(\tau )} \right)d\tau } } \right)} \right. \\
&\left. { - \frac{1}{{\Gamma \left( \alpha  \right)}}\int\limits_a^x {\psi '\left( t \right)} \left( {\psi \left( x \right) - \psi \left( t \right)} \right)^{\alpha  - 1} f\left( {t,0,\int\limits_a^t {w\left( {t,\tau ,0} \right)d\tau } } \right)} \right| \\
& \le p_2  + \frac{1}{{\Gamma \left( \alpha  \right)}}\int\limits_a^x {\psi '\left( t \right)} \left( {\psi \left( x \right) - \psi \left( t \right)} \right)^{\alpha  - 1}  \\ 
& \times \left| {f\left( {t,\mathfrak{z}(t),\int\limits_a^t {w\left( {t,\tau ,\mathfrak{z}(\tau )} \right)d\tau } } \right) - f\left( {t,0,\int\limits_a^t {w\left( {t,\tau ,0} \right)d\tau } } \right)} \right| \\
& \le p_2  + \frac{1}{{\Gamma \left( \alpha  \right)}}\int\limits_a^x {\psi '\left( t \right)} \left( {\psi \left( x \right) - \psi \left( t \right)} \right)^{\alpha  - 1}  \\ 
&\times \left[ {Q_3 \left( t \right)\left( {\left| {\mathfrak{z}\left( t \right)} \right| + \int\limits_a^t {Q_4 \left( \tau  \right)\left| {\mathfrak{z}\left( \tau  \right)} \right|d\tau } } \right)} \right]dt.
\tag{4.5}
\end{align*}
Applying the inequality in Lemma $2.3$ to above equation $(4.5)$ with $\mathfrak{u}(t) = \mathfrak{u}(x)$, $ \mathfrak{u}_0  = p_2 $, $f\left( s \right) = \frac{{\psi '\left( t \right)\left( {\psi \left( x \right) - \psi \left( t \right)} \right)^{\alpha  - 1} }}{{\Gamma \left( \alpha  \right)}}Q_3 \left( t \right) $ and $g\left( \tau  \right) = Q_4 \left( t \right) $  we get the required inequality $(4.4)$.
\paragraph{\textbf{Remark}} The estimate obtained in above $(4.4)$ gives the bound for the solution of eqution $(1.1)$. If $(4.4)$ is bounded then $(1.1)$ is also bounded.

\section{Continuous dependence of  solution}
Now here we present the result on continuous dependence property of solution of equation $(1.1)-(1.2)$ on the function involved therein.
\paragraph{} Consider the equation $(1.1)-(1.2)$ and the corresponding equation
\[
{}^HD_{a + }^{\alpha ,\beta ;\psi } \mathfrak{z}(x) = \overline f \left( {x,\mathfrak{z}\left( x \right),\int\limits_a^x {\overline w \left( {x,t,\mathfrak{z}\left( t \right)} \right)dt} } \right),
\tag{5.1}\]
\[
I_{a + }^{1 - \gamma ;\psi } \mathfrak{z}\left( a \right) = \mathfrak{z}_a ,\,\,\,\,\,\,\,\,\,\,\,\,\,\,\,\gamma  = \alpha  + \beta \left( {1 - \alpha } \right),
\tag{5.2}\]
where $\overline f :\left[ {a,b} \right] \times \mathbb{R} \times \mathbb{R} \to \mathbb{R}$ and $\mathfrak{z}_a$ is constant and $\overline w :\left[ {a,b} \right] \times \left[ {a,b} \right] \times \mathbb{R} \to \mathbb{R}$.
\paragraph{} Now we present the theorem which gives the continuous dependence of solution of $(1.1)$.

\paragraph{\textbf{Theorem 5.1.}} Suppose that the functions $f$ and $w$ in $(1.1)$ and $(5.1)$ satisfy the conditions $(4.1)$ and $(4.2)$. Let $v(x)$ be a solution of $(5.1)$ and 
\begin{align*}
&\left| {\mathfrak{z}_a  - \mathfrak{v}_a } \right|\frac{{\left| {\left( {\psi \left( x \right) - \psi \left( a \right)} \right)^{\gamma  - 1} } \right|}}{{\Gamma \left( \gamma  \right)}} + \frac{1}{{\Gamma \left( \gamma  \right)}}\int\limits_a^x {\psi '\left( t \right)} \left( {\psi \left( x \right) - \psi \left( a \right)} \right)^{\alpha  - 1}  \\ 
&\left| {f\left( {t,\mathfrak{z}(t),\int\limits_a^t {w\left( {t,\tau ,\mathfrak{z}\left( \tau  \right)} \right)d\tau } } \right) - \overline f \left( {t,\mathfrak{v}(t),\int\limits_a^t {\overline w \left( {t,\tau ,\mathfrak{v}\left( \tau  \right)} \right)d\tau } } \right)} \right| < \varepsilon ,
\tag{5.3}
\end{align*}
then
\begin{align*}
\left| {\mathfrak{z}(x) - \mathfrak{v}(x)} \right| 
& \le \varepsilon \left[ {1 + \int\limits_a^x {\psi '\left( t \right)} \left( {\psi \left( x \right) - \psi \left( a \right)} \right)^{\alpha  - 1} Q_3 } \right.\left( t \right) \\ 
&\left. {\exp \left( {\int\limits_a^t {\left( {\frac{{\psi '\left( s \right)\left( {\psi \left( x \right) - \psi \left( s \right)} \right)^{\alpha  - 1} }}{{\Gamma \left( \gamma  \right)}}Q_3 \left( s \right) + Q_4 \left( s \right)} \right)ds} } \right)dt} \right], 
\tag{5.4}
\end{align*}
where $f,w$ and $\overline f, \overline w$ are the functions involved in equation $(1.1)$ and $(5.1)$ and $\epsilon >0$ is arbitrarily small constant. Then the solution $\mathfrak{z}(x)$ of equation $(1.1)$ depends continuous on functions involved in $(1.1)$.
\paragraph{\textbf{Proof.}}
Let $\mathfrak{u}(x)=|\mathfrak{z}(x)-\mathfrak{v}(x)|$, $x \in [a,b]$. Since $\mathfrak{z}(x)$ and $\mathfrak{v}(x)$ are solutions of equation $(1.1)$ and $(5.1)$ and from the hypotheses we get

\begin{align*}
\left| {\mathfrak{u}\left( x \right)} \right|
&\le \left| {\mathfrak{z} - \mathfrak{v}} \right|\frac{{\left| {\left( {\psi \left( x \right) - \psi \left( a \right)} \right)^{\gamma  - 1} } \right|}}{{\Gamma \left( \gamma  \right)}} + \frac{1}{{\Gamma \left( \gamma  \right)}}\int\limits_a^x {\psi '\left( t \right)} \left( {\psi \left( x \right) - \psi \left( a \right)} \right)^{\alpha  - 1}  \\
&\times \left| {f\left( {t,\mathfrak{z}\left( t \right),\int\limits_a^t {w\left( {t,\tau ,\mathfrak{z}\left( \tau  \right)} \right)d\tau } } \right) - f\left( {t,\mathfrak{v}(t),\int\limits_a^t {w\left( {t,\tau ,\mathfrak{v}\left( \tau  \right)} \right)d\tau } } \right)} \right| \\ 
&+ \frac{1}{{\Gamma \left( \gamma  \right)}}\int\limits_a^x {\psi '\left( t \right)} \left( {\psi \left( x \right) - \psi \left( a \right)} \right)^{\alpha  - 1}  \\ 
& \times \left| {f\left( {t,\mathfrak{v}\left( t \right),\int\limits_a^t {w\left( {t,\tau ,\mathfrak{v}\left( \tau  \right)} \right)d\tau } } \right) - \overline f \left( {t,\mathfrak{v}(t),\int\limits_a^t {\overline w \left( {t,\tau ,\mathfrak{v}\left( \tau  \right)} \right)d\tau } } \right)} \right| \\ 
&\le \varepsilon  + \frac{1}{{\Gamma \left( \gamma  \right)}}\int\limits_a^x {\psi '\left( t \right)} \left( {\psi \left( x \right) - \psi \left( a \right)} \right)^{\alpha  - 1}  \\ 
&\times \left[ {Q_3 \left( t \right)\left| {\mathfrak{z}\left( t \right) - \mathfrak{v}\left( t \right)} \right| + \int\limits_a^t {Q_4 \left( s \right)\left| {\mathfrak{z}\left( s \right) - \mathfrak{v}\left( s \right)} \right|ds} } \right]. 
\tag{5.5}
\end{align*}
Now an application of lemma $(2.3)$ to $(5.5)$ and similar as in proof of Theorem $(4.1)$  we get the required inequality $(5.3)$.


\begin{thebibliography}{999}
\bibitem {Abb}   S. Abbas, M. Benchohra ,  J. Lazreg and Y. Zhou, 
\newblock  A survey on Hadamard and Hilfer fractional differential equations: Analysis and stability,  
\newblock  {\em  Chaos Solitons  Fract. },  102 (2017), 47–71.


\bibitem {Alm1} R. Almeida,
\newblock  A Caputo fractional derivative of a function with respect to another function,
\newblock {\em Commun. Nonlinear Sci. Numer. Simult.}, 44(2017), 460-481.

\bibitem {Fur}  K.M. Furati, M.D. Kassim and N.e-. Tatar,
\newblock    Existence and uniqueness for a problem involving Hilfer fractional derivative,
\newblock  {\em Comput.  Math. Appl. },  64 (2012), 1616–1626.


\bibitem {Bal} D. Baleanu and A. M. Lopes,   
\newblock Handbook of Fractional Calculus with Applications, Applications in Engineering, Life and Social Sciences, 
\newblock {\em De Gruyter.}, Vol. 7:2019.

\bibitem {Har1} S. Harikrishnan, K. Shah  and K. Kanagarajan ,
\newblock    Study of a boundary value problem for fractional order $\psi$-Hilfer fractional derivative   ,
\newblock  {\em  Arab. J. Math.}, (2019),  .

\bibitem {Har2} S. Harikrishnan , K. Kanagarajan and D. Vivek
\newblock   Existence and stability results for boundary value problem for differential equation with  $\psi$-Hilfer fractional derivative  ,
\newblock  {\em Journal of Applied Nonlinear Dynamics },  8(2),(2019), 251-259.

\bibitem {Kar} K. Karthikeyan and J.J. Trujillo,
\newblock    Existence and uniqueness results for fractional integrodifferential
equations with boundary value conditions,
\newblock  {\em Commun. Nonlinear Sci. Numer. Simulat.}, 17 (2012), 4037–4043.

\bibitem {Kil1} A. A. Kilbas, H. M. Srivastava and J. J. Trujilio,   
\newblock Theory and applications of fractional differential equations,
\newblock {\em North Holland Mathematics studies.}, Vol 207, 2006.




\bibitem {Van2} K. D. Kucche, A. D. Mali and J. C. Sousa,
\newblock On the nonlinear $\psi$-Hilfer fractional differential equations.
\newblock {\em Comp. Appl. Math.},  (2018), 38-73.

\bibitem {Bgp1} B.G. Pachpatte,
\newblock  A Note of Gronwall-Bellman Inequality,
\newblock {\em J Math. Anal. Appl.}, 44(1973), 758-762.

\bibitem {Pet} I. Petras,   
\newblock Handbook of Fractional Calculus with Applications, Applications in Control, 
\newblock {\em De Gruyter.}, Vol. 6:2019.


\bibitem {Sam1} S.G. Samko, A. A. Kilbas and O. I. Marichev, 
\newblock Fractional Integrals and Derivatives: Theory and Applications,
\newblock {\em  Gordon and Breach Science Publishers}, 1993. 


\bibitem {Van1} J. C. Sousa and E. C. De Oliveira,
\newblock On the $\psi$-Hilfer fractional derivative,
\newblock {\em Commun Nonlinear Sci. Numer. Simultat}, 60(2018), 72-91.


\bibitem {Van3} J. C. Sousa, F.G. Rodrigues and E. C. De Oliveira,
\newblock Stability of the fractional Volterra integrodifferential equation by means of $\psi$-Hilfer Operator,
\newblock {\em Math. Meth. Appl. Sci.}, 42(2019), 3033-3043.

\bibitem {Van4} J. C. Sousa and E. C. De Oliveira,
\newblock On the stability of a Hyperbolic Fractional Partial Differential Equation,
\newblock {\em Differ. Equ. Dyn. Syst.}, (2018). 

\bibitem {Van5} J. C. Sousa and E. C. De Oliveira,
\newblock A Gronwall inequality and the Cauchy type problem by means of $\psi$-Hilfer Operator,
\newblock {\em Differ. Equ. Appl.}, 11(1), (2019), 87-106.

\bibitem {Sub}  R. Subashini ,  C. Ravichandran,  K. Jothimani and H. M. Baskonus,
\newblock    Existence results of Hilfer integrodifferential equations with fractional order,
\newblock  {\em Discrete  Cont. Dyn-S }, 13(3)(2020), 911-923.

\bibitem {Tat} S. Tate, V. V. Kharat and H. T. Dinde
\newblock  A nonlinear mixed fractional integrodifferential equation with positive constant coefficient,
\newblock  {\em  Filomat}, 33:17,  (2019), 5623-5638.  .

\bibitem {Tha}S. Thabet, B.  Ahmad and  R. P. Agrwal,
\newblock    On abstract Hilfer fractional integrodifferential equations withboundary conditions,
\newblock  {\em Arab J. Math. Sci. },  (2019).

\bibitem {Wan}  J. Wang and Y. Zhang,
\newblock     Nonlocal initial value problems for differential equations with Hilfer fractional derivative,
\newblock  {\em  Appl. Math. Comput. }, 266(2015), 850–859.

\bibitem {Zha} L. Zhang, B. Ahmad, G. Wang and R. P. Agrwal,
\newblock  Nonlinear fractional integrodifferential equations on unbounded domains in a Banach space,
\newblock  {\em  J. Comput.  Appl. Math.},  249(2013), 51-56.


\end{thebibliography}
\end{document}